\documentclass[11pt]{amsart}
\usepackage{amsfonts}
\usepackage[all]{xy}
\begin{document}

\newcommand{\bk}{{\bf k}}
\newcommand{\bl}{{\bf l}}
\newcommand{\bn}{{\bf n}}
\newcommand{\bm}{{\bf m}}
\newcommand{\C}{{\mathbb C}}
\newcommand{\ok}{{\overline{k}}}
\newcommand{\e}{{\bf e}}
\newcommand{\he}{{\hat{\bf e}}}
\newcommand{\SR}{{\rm SR}}
\newcommand{\tSR}{{\widetilde{\rm SR}}}
\newcommand{\NSymm}{{\sf NSymm}}
\newcommand{\QSymm}{{\sf QSymm}}
\newcommand{\Z}{{\mathbb Z}}
\newcommand{\jd}{{\rm DJ}}
\newcommand{\N}{{\mathbb N}}
\newcommand{\bL}{{\overline{\Lambda}}}
\newcommand{\bTh}{{\overline{\Theta}}}
\newcommand{\bV}{{\overline{V}}}
\newcommand{\la}{{\langle}}
\newcommand{\ra}{{\rangle}}
\newcommand{\T}{{\mathbb T}}
\newcommand{\U}{{\rm U}}
\newcommand{\Q}{{\mathbb Q}}
\newcommand{\res}{{\rm res}}
\newcommand{\Spec}{{\rm Spec \;}}
\newcommand{\Hom}{{\rm Hom}}

\title{The Baker-Richter spectrum as cobordism of quasitoric manifolds}
\author{Jack Morava}
\address{Department of Mathematics, Johns Hopkins University, Baltimore,Maryland}
\email{jack@math.jhu.edu}
\author{Nitu Kitchloo}
\address{Department of Mathematics, Johns Hopkins University, Baltimore,Maryland}
\email{nitu@math.jhu.edu}

\date {14 January 2012}

\begin{abstract}{Baker and Richter construct a remarkable $A_\infty$ ring-spectrum $M\Xi$ 
whose elements possess characteristic numbers associated to quasisymmetric functions; its 
relations, on one hand to the theory of noncommutative formal groups, and on the other to 
the theory of omnioriented (quasi)toric manifolds [in the sense of Buchstaber, Panov, and 
Ray], seem worth investigating.}\end{abstract}

\maketitle
\bigskip

\noindent
{\bf Introduction:} Most of this paper is a draft for a talk JM wishes he had
given at the August 2011 conference on toric manifolds [10] at Queen's University,
Belfast (as opposed to the talk he actually gave). Thanks to Thomas H\"uttemann 
for organizing that very interesting meeting, and to Tony Bahri, Martin Bendersky, 
Fred Cohen, and Sam Gitler for helpful conversations there. Those notes are little 
more than a collage of conversations, suggestions, and howler-preventing interventions 
courtesy of Andy Baker, Michiel Hazewinkel, Birgit Richter, and Taras Panov, in the 
course of the last few years; he is deeply indebted to them all. \bigskip

\noindent
He is also extremely grateful to {\bf Nitu Kitchloo}, for permission to include 
an appendix by the latter, which outlines some work in progress. \bigskip

\section{Conjectures about the spectrum $M\Xi$} \bigskip

\noindent
{\bf 1.1} This is the Thom spectrum defined by A Baker and B Richter [2, with
slightly different notation] constructed by pulling back the canonical bundle $\xi \to
B\U$ along the abelianization map
\[
\Omega \Sigma B\T = B \Omega^2 \Sigma B\T \to B\U \;.
\]
[An analytic construction of the associated representation
\[
\Xi := \Omega^2 \Sigma B\T \to \U
\]
might be very interesting.] \bigskip

\noindent
They show that $M\Xi$ is an $A_\infty$ ring-spectrum, with $M\Xi_*$ torsion-free and 
concentrated in even degrees, and that the Hurewicz homomorphism
\[
M\Xi_* \to H_*(\Omega \Sigma \C P^\infty) \cong \NSymm_*
\]
takes values in the (graded) ring of noncommutative symmetric functions [6 \S 4.2, 8];
it is injective, and becomes an isomorphism after tensoring with $\Q$. Finally, and most 
striking of all, they show that $M\Xi \otimes Z_{p}$ is a wedge of 
copies of $BP$. \bigskip

\noindent
{\bf 1.2} A (unital) $S$-algebra $A$ defines a cosimplicial algebra
\[
A^\bullet : n \to A^{\wedge n} 
\]
with maps built from its unit and multiplication. In good cases [1] this is a resolution, 
in a suitable sense, of the sphere spectrum, and its homotopy groups define a cosimplicial algebra
\[
\xymatrix{
A_* \ar@<-.5ex>[r] \ar@<.5ex>[r] & A_*A \otimes_A A_*A \ar[r] \ar@<1ex>[r] \ar@<-1ex>[r] & \dots}
\]
which leads to the construction of an Adams spectral sequence. \bigskip

\noindent
If $A = M\U$, then
\[
  \pi_*(M\U \wedge M\U) = M\U_* \otimes_\Z S_*
\]
is the product of the Lazard ring with the algebra $S_* \cong H_*(B\U,\Z)$ of functions on the group
of formal power series under composition, and the resulting cosimplicial ring can be interpreted as
a presentation
\[
MU_* \to M\U_* \otimes S_* \to (M\U_* \otimes S_*) \otimes_{M\U} (M\U_* \otimes S_*)
\to \dots
\]
of (the graded algebra of functions on) the moduli stack of one-dimensional formal group laws. The 
classical Steenrod augmentation $M\U_* \to \Z$ classifies the additive group law, and its 
composition
\[
M\U_* \to M\U_* \otimes S_* \to \Z \otimes S_* = S_*
\]
with the coaction represents inclusion
\[
\Spec S_* \to \Spec M\U_*
\]
of the orbit of the additive group, under coordinate changes, in the moduli stack of formal
groups. It is also the (injective) Hurewicz map
\[
\pi_*(M\U) \to H_*(M\U,\Z) \cong \Hom(H^*(B\U),\Z)
\]
which assigns to a manifold, its collection of Chern numbers; and from either point of view 
it is a rational isomorphism. \bigskip

\noindent
The cosimplicial algebra $\pi_*M\U^\bullet \otimes \Q$ is thus a resolution of $\pi^S_* \otimes 
\Q = \Q$; it is a cosimplicial presentation of the stack over $\Q$ defined by 
the action of the group of formal diffeomorphisms on itself\begin{footnote}{If $G$ is a group 
object, then the category $[G/G]$ defined by its translation action is equivalent to the category 
with one object and its identity morphism.}\end{footnote}. \bigskip

\noindent
{\bf 1.3} It would be nice to have a similar description for $A = M\Xi$, when
\[
\pi_*(M\Xi \wedge M\Xi) \; = \; M\Xi_*M\Xi \; \cong \; M\Xi_* \otimes_\Z \NSymm_* \;,
\]
but we don't yet have a good description of the coaction maps in the cosimplicial ring it
defines. I am indebted to Michiel Hazewinkel for suggesting the following possibility:
\bigskip

\noindent
{\bf Conjecture:} The Hurewicz homomorphism
\[
M\Xi_*M\Xi \to H_*(M\Xi \wedge M\Xi) \cong \NSymm_* \otimes_\Z \NSymm_*
\]
is a homomorphism of Hopf algebras, with target the Novikov 
double [9] of the Hopf algebra defined by the diagonal
\[
\Delta_{BFK} Z(t) = \res_{u=0} \; Z(u) \otimes (u - Z(t))^{-1}
\]
[4 \S 2.4] on the ring
\[
\NSymm_* = \Z \la Z_i \:|\: i \geq 1 \ra, \; Z(t) = \sum_{i \geq 0} Z_i t^{i+1}
\]
of noncommutative symmetric functions. \bigskip

\noindent
{\bf 1.4} In further work Baker and Richter construct [3] an injective homomorphism $\lambda_{BR}$
\[
z \mapsto c + \sum_{i > 0} z_i c^{i+1} : M\Xi^*B\T \cong M\Xi_*\la \la z \ra \ra \to H_*M\Xi
[[c]] \cong \NSymm_* [[c]]
\]
of Hopf algebras, where $c$ is a central element corresponding to the Chern class for line
bundles in ordinary cohomology. The diagram
\[
\xymatrix{
M\Xi_* \ar[d] \ar[r] & \NSymm_* \ar@{.>}[dl] \ar[d] \\
M\Xi_* \otimes \Q \ar[r]^{\lambda_{BR} \otimes \Q} & \NSymm_* \otimes \Q }
\]
lets us regard the coefficients $z_i$ of their logarithm as elements of $M\Xi_{2i}
\otimes \Q$.\bigskip

\noindent
{\bf 1.5 Conjecture:} With the left vertical homomorphism defined by the natural coaction, the diagram
\[
\xymatrix{
M\Xi^*B\T \ar[r]^{\lambda_{BR}} \ar[d] & \NSymm_*[[c]] \ar[d]^{\Delta_{BFK}} \\
M\Xi^*B\T \hat{\otimes}_{M\Xi_*} M\Xi_*M\Xi \cong M\Xi^*B\T \hat{\otimes}_\Z \NSymm_* \ar[r]^{\Delta_{B\T}} 
& (\NSymm_* \otimes_\Z \NSymm_*)[[c]] }
\]
commutes; where 
\[
\Delta_{B\T}(1 \otimes c) := 1 \otimes Z(c) : \NSymm_*[[c]] \to (\NSymm_* \otimes_\Z \NSymm_*)[[c]] \;.
\]
Moreover, $\lambda_{BR} \otimes \Q$ maps $\pi_*M\Xi^\bullet \otimes \Q$ isomorphically to the resolution
\[
\xymatrix{
\Q \ar[r] & \NSymm_* \otimes \Q \ar@<-.5ex>[r] \ar@<.5ex>[r] & (\NSymm_* \otimes_\Q \NSymm_*) \ar[r] \ar@<1ex>[r] 
\ar@<-1ex>[r] & \dots }
\]
\bigskip

\section{Characteristic numbers for quasitoric manifolds}\bigskip

\noindent
{\bf 2.1} A complex-oriented $2m$-dimensional quasitoric manifold $M$ has a
simple quotient polytope $P$, with an omniorientation [5 \S 5.31]
defined by a characteristic map $\Lambda$ [5 \S 5.10] from the ordered set
of vertices of the simplicial $(m-1)$-sphere $K_P$ bounding the
dual simplicial complex $P^*$ [5 \S 1.10], to a free abelian group $\Theta$
with generators $\theta_1, \dots, \theta_m$. \bigskip

\noindent
The cohomology of $M$ can be naturally identified [5 \S 5.2.2, \S 6.5, 7] with 
a quotient
\[
k^*(K_P) \otimes_{P(\Theta^*)} \Z \;,
\]
of the Stanley-Reisner face ring [5 \S 3.1, 3.4]; it is an algebra over
the symmetric algebra $P(\Theta^*)$ on the $\Z$-dual of $\Theta$ via 
\[
\Lambda^* : P(\Theta^*) \to P(V_P^*) \to k^*(K_P) 
\]
(where $V_P$ is the free abelian group generated by the vertices of $K_P$,
and $\Z$ is a $P(\Theta^*)$-algebra via augmentation). \bigskip

\noindent
{\bf 2.2} The order on the vertex set of $K_P$ embeds it as the initial segment
of the natural numbers $\N$, identifying $V_P$ with a subgroup of a free
abelian group $\bV$ on a countable set of generators. A similar identification
embeds $\Theta$ in another free abelian group $\bTh$ on a countable set of 
generators, defining an extension
\[
\bL : \bV \to \bTh
\]
of $\Lambda$ by $\bL([i]) = \theta_i$ when $i$ is {\bf not} a vertex of $K_P$. 
The resulting homomorphism
\[
P(\bTh^*) \to P(\bV^*) = P(V^*_P) \otimes P(\bV^*_{\neg P})
\]
makes
\[
\ok^*(K_P) := k^*(K_P) \otimes_\Z P(\bV_{\neg P}^*) 
\]
into an algebra over the polynomial ring $P(\bTh^*)$ generated by a 
countable sequence of variables, such that
\[
H^{2*}(M,\Z) \cong \ok^*(K_P) \otimes_{P(\bTh^*)} \Z \;.
\] \bigskip

\noindent
{\bf 2.3} For an {\bf ordered} partition 
\[
\bm \; = \; m_1 \; +  \cdots + \;  m_r
\]
of $m$, let $[\bm]_P$ be the image in $H^{2m}(M,\Z)$ of the formal sum 
\[
[\bm] := \sum_{i_1 < \dots < i_r} x_{i_1}^{m_1} \cdots x_{m_r}^{n_r} \;, 
\]
where $x_i$ is the polynomial generator corresponding to $i \in \bV$ [6 \S 4]. If $i$ is 
sufficiently small, $x_i$ corresponds to a vertex of $K_P$; otherwise, it is a kind
of dummy element, and is killed by $-\otimes_{P(\bTh^*)} \Z$. \bigskip

\noindent
{\bf 2.4} More generally, if $M$ and $N$ are almost-complex quasitoric manifolds
of dimension $m,n$ respectively, with quotient polytopes $P,Q$, then the product
\[
k^*(K_P) \otimes_\Z k^*(K_Q) \cong k^*(K_P * K_Q)
\]
of face rings is naturally isomorphic to the face ring of the join 
\[
K_P * K_Q \cong K_{P \times Q}
\]
of the simplicial spheres $K_P$ and $K_Q$ [5 \S 2.13]. \bigskip

\noindent
{\bf Claim:} the corresponding isomorphism
\[
H^*(M) \otimes H^*(N) \to H^*(M \times N)
\]
sends $[\bm]_P \otimes [\bn]_Q$ to 
\[
[\bm + \bn]_{P \times Q} := [m_1 + \dots + m_r + n_1 + \dots + n_s]_{P \times Q} \;.
\] \bigskip

\noindent
{\bf Proof:} $[\bm + \bn]_{P \times Q}$ is the image in $H^{2(m+n)}(M \times N)$ of 
\[
\sum_{i_1 < \dots < i_{r+s}} x_{i_1}^{m_1} \cdots x_{i_{r+s}}^{n_s} \;,
\]
summed over strings $i_1 < \dots < i_{r+s}$ of elements of the disjoint union of the 
vertex sets of $P^*$ and $Q^*$. For a monomial of this sort to have a nontrivial image 
in the top-dimensional cohomology of $M \times N$, the elements of the set $\{ x_{i_1}, 
\dots, x_{i_r} \}$ must be vertices of $K_P$, and those of $\{ x_{i_{r+1}}, \dots, 
x_{i_{r+s}} \}$ must be vertices of $K_Q$. The image of the sum is thus the product of 
the images of the sums $[\bm]_P$ and $[\bn]_Q$ (modulo the identification of the 
top-dimensional cohomology group of $P$ (resp. $Q$) with the integers). $\Box$ \bigskip

\noindent
This construction associates to a $2m$-dimensional complex-oriented quasitoric 
manifold $M$, a homomorphism
\[
\bm \mapsto [\bm]_P : \QSymm^m \to \Z 
\]
of abelian groups (ie a noncommutative symmetric function ${\bf M}$), which sends $M \times N$ 
to a noncommutative symmetric function
\[
{\bf M \times N} \; = \; {\bf M} \bullet {\bf N}
\]
equal to the product of the noncommutative symmetric functions {\bf M} and {\bf N}. \bigskip

\noindent
{\bf Problem:} What is the noncommutative symmetric function $\bf{CP}_n$ defined by complex
projective $n$-space, with its usual toric structure (and the $n$-simplex as associated
polytope)? \bigskip

\noindent
{\bf 2.5} This might be paraphrased as saying that the Davis-Januszkiewicz construction defines 
a ring homomorphism from the algebra generated by the monoid, under join, of certain
omnioriented simplicial spheres, to the free graded associative algebra $\NSymm_*$; in 
other words, something like a coordinate patch for a noncommutative space 
of quasitoric manifolds. \bigskip

\section{Appendix, by Nitu Kitchloo:} \bigskip

\noindent
A quasitoric manifold $M$ of dimension $2m$ admitting an action of a torus $T$ of rank $m$
is associated with a polytope $P$. Assume $F = \{ f_i \}$ is the set consisting of the co-dimension 
one faces $f_i$ of $P$. The data required to construct $M$ involves a collection of primitive 
characteristic weights $\lambda_i \in \pi_1(T)$, indexed on the set $F$. \bigskip

\noindent
Let $\hat{T}$ denote the torus $(S^1)^F$, of rank given by the cardinality of $F$, with a 
canonical set of generating circles indexed by the faces $f_i$. Let $H \subset \hat{T}$ 
denote the kernel of the map $\lambda : \hat{T} \longrightarrow T$, defined by $\lambda(\exp(tf_i)) 
= \exp(t\lambda_i)$. \bigskip

\noindent
The procedure for constructing $M$ is described as follows: Notice that $\hat{T}$ acts on $\C^F$ 
in a canonical way via Hamiltonian symplectomorphisms. This induces an action of $H$. Let $\mathcal{H}$ denote 
the Lie algebra of $H$. Let $\varphi : \C^F \longrightarrow \mathcal{H}^\ast$ denote the moment 
map of the $H$ action. Let the "Moment angle complex" $Z(P)$ denote the preimage of a regular value. The manifold $M$ is defined as the orbit space: $M = Z(P)/H$. From this it follows easily that: \bigskip

\noindent
{\bf 3.1 Claim:} $M = Z(P)/H$ has stable tangent bundle classified by the composite map:
\[ 
\tau(M) : Z(P)/H \longrightarrow BH \longrightarrow B\hat{T} \longrightarrow BU(F), 
\]
where $BU(F)$ denotes the group of unitary transformations of $\C^F$, with maximal torus $\hat{T}$.
\bigskip

\noindent
Let us now try to find a natural Thom spectrum that is the receptacle for the cobordism class of $M$. 
Firstly notice that there is a commutative diagram:
\[
\xymatrix{
B\hat{T} \ar[d]   \ar[r] & \Omega \Sigma BU(1)  \ar[d]^{\varphi} \\
BU(F) \ar[r] & BU
}
\]
where $\varphi : \Omega \Sigma BU(1) \longrightarrow BU$ is the $A_{\infty}$ extension of the 
inclusion map $BU(1) \longrightarrow BU$, and the map $B\hat{T} \longrightarrow \Omega \Sigma 
BU(1) $ is the inclusion of the $|F|$-th James filtration. We conclude:\bigskip

\noindent
{\bf 3.2 Corollary:} Let $\overline{M}\Xi$ denote the Thom spectrum of $-\varphi$, then the cobordism 
class of $M$ belongs to $\pi_{2m} \overline{M}\Xi$.
\bigskip

\noindent 
{\bf 3.3 Remark:} Notice that since $\varphi$ is an $A_{\infty}$ map, the Thom spectrum of $-\varphi$ 
is equivalent to the Thom spectrum $M\Xi$ of $\varphi$, as seen easily from the following commutative diagram:
\[
\xymatrix{
\Omega \Sigma BU(1)  \ar[d]^{\varphi}  \ar[r]^{-Id} & \Omega \Sigma BU(1)  \ar[d]^{-\varphi} \\
BU \ar[r]^{Id} & BU
}
\]
In particular, the algebraic procedure described in Section 2 is indeed a method of computing the 
characteristic numbers of the tangent bundle of $M$. 
\bigskip

\bibliographystyle{amsplain}

\end{document}